\documentclass{amsart}
\usepackage{amsmath, amsthm, amssymb, amscd}

\DeclareMathOperator{\as}{asdim} 
 
 \DeclareMathOperator{\dist}{dist}
 
\DeclareMathOperator{\supp}{supp}

\begin{document}
\author{Gregory C. Bell}
\date{November 26, 2002}
%\commby{Stephen D. Smith}

\address{Department of Mathematics, 328 Natural Sciences Building,
Louisville, KY 40292}

\email{gcbell01@gwise.louisville.edu}

\title[Asymptotic properties of groups acting on complexes]{Asymptotic properties of groups acting on complexes}

\newcommand{\To}{\longrightarrow}
\newcommand{\sA}{\mathcal{A}}
\newcommand{\sB}{\mathcal{B}}
\newcommand{\sU}{\mathcal{U}}
\newcommand{\sV}{\mathcal{V}}
\newcommand{\sW}{\mathcal{W}}
\newcommand{\sx}{\mathcal{X}}
\newcommand{\sY}{\mathcal{Y}}
\newcommand{\sZ}{\mathcal{Z}}
\newcommand{\Real}{\mathbf{R}}
\newcommand{\Complex}{\mathbf{C}}
\newcommand{\so}{\Rightarrow}
\newcommand{\sX}{\mathcal{X}}
%\newcommand{\as}{\textrm{asdim}}
% THEOREM ENV. ---------------------------------------------------
\newtheorem{Theorem}{Theorem}
\newtheorem*{thm}{Theorem}
\newtheorem*{cor}{Corollary}
\newtheorem{Lemma}[Theorem]{Lemma}
\newtheorem{Proposition}[Theorem]{Proposition}
\newtheorem*{defn}{Definition}

\begin{abstract} We examine asymptotic dimension and property A
for groups acting on complexes.  In particular, we prove that the
fundamental group of a finite, developable complex of groups will
have finite asymptotic dimension provided the geometric
realization of the development has finite asymptotic dimension and
the vertex groups are finitely generated and have finite
asymptotic dimension.  We also prove that property A is preserved
by this construction provided the geometric realization of the
development has finite asymptotic dimension and the vertex groups
all have property A. These results naturally extend the
corresponding results on preservation of these large-scale
properties for fundamental groups of graphs of groups. We also use
an example to show that the requirement that the development have
finite asymptotic dimension cannot be relaxed.
\end{abstract}

\maketitle

\section{Introduction}

The asymptotic approach to finitely generated groups became
popular following the work of Gromov \cite{Gr1}.  In his study of
asymptotic invariants of finitely generated groups, Gromov defined
asymptotic dimension ($\as$), the large-scale analog of Lebesgue
covering dimension. G. Yu \cite{Yu1} applied asymptotic dimension
to the Novikov higher signature conjecture for groups, showing
that the conjecture holds for groups with finite asymptotic
dimension. Later, Yu \cite{Yu2} defined another asymptotic
invariant for discrete metric spaces and finitely generated groups
called {\em property A.} This is a weak form of amenability which
also implies the Novikov conjecture for groups.  (For an
introduction to the Novikov and related conjectures, see
\cite{FRR}.)

We wish to consider finitely generated groups as metric spaces.
Let $\Gamma$ be a finitely generated group with generating set
$S=S^{-1}.$ The {\em $S$-norm} on $\Gamma$ is the norm given by
setting $\|\gamma\|_S=0$ precisely when $\gamma$ is the group
identity and otherwise taking $\|\gamma\|_S$ to be the minimal
length of any $S$-word presenting the element $\gamma.$  Then, one
can define the (left-invariant) {\em word metric} associated to
$S$ by $\dist_S(g,h)=\|g^{-1}h\|_S.$  The metrics corresponding to
two finite generating sets $S$ and $S'$ are Lipschitz equivalent.
Asymptotic dimension and property A are invariants of Lipschitz
equivalent metric spaces, so these properties are intrinsic to the
group $\Gamma$ and not the metric space associated to a specific
generating set.

In view of Yu's results, \cite{Yu1},\cite{Yu2}, it is important to
know which groups have finite $\as$ or property A.  Gromov
\cite{Gr1} showed that hyperbolic groups have finite $\as.$
Dranishnikov and Januszkiewicz proved in \cite{DJ} that Coxeter
groups have finite $\as.$ Dranishnikov and the author proved that
the finiteness of $\as$ is preserved by the amalgamated free
product and HNN extension, and more generally any fundamental
group of a finite graph of groups with vertex groups having finite
$\as,$ see \cite{BD1}, \cite{BD2}.  Higson and Roe \cite{HR}
showed that finitely generated groups with finite asymptotic
dimension have property A.  In \cite{Tu}, J.-L. Tu proved that
property $A$ is preserved by the fundamental group of a finite
graph of groups where the vertex groups all have property A.  At
first it was not know whether there could be a finitely generated
group with infinite $\as$ or which does not have property A. A
recent example of such a group due to Gromov \cite{Gr2} and
Dranishnikov, Gong, Lafforgue and Yu \cite{Detal} , have made
determining precisely which groups have these properties
interesting.

The Bass-Serre theory of graphs of groups generalizes the
constructions of amalgamated free products and HNN extensions (see
\cite{S}).  There is a direct correspondence between groups acting
without inversion on trees and fundamental groups of graphs of
groups. Complexes of groups were introduced by Haefliger
\cite{Hae} in order to describe actions of groups on simply
connected simplicial complexes in the same way that graphs of
groups describe the action of groups on trees. The problem that
arises is that the quotient of a simplicial complex by a
simplicial action may identify faces of simplices.  So it may not
be the case that the quotient is a simplicial complex. This
problem is avoided by introducing combinatorial substitutes for
simplicial complexes called small categories without loops
(scwols).

In the second section we develop the necessary theory of complexes
of groups following \cite{BH}. We define scwols, group actions on
scwols, complexes of groups, developability of complexes of
groups, and the associated fundamental group of a complex of
groups.

In the third section, we define the $R$-stabilizer which plays the
r\^{o}le of the stabilizer in the study of the large-scale
properties of groups acting on metric spaces.  We also
characterize the $R$-stabilizers for groups acting on scwols.

In the fourth section, we define $\as$ and obtain our main result
on asymptotic dimension:

\begin{thm} Let $G(\sY)$ be a developable complex of groups over a
finite scwol $\sY$ with development $\sX$ such that
$\as\sX<\infty,$ and such that every local group $G_\sigma$ has
$\as G_\sigma<\infty.$  If $\pi_1$ denotes the fundamental group
of the complex of groups, we have $\as \pi_1<\infty.$
\end{thm}

In the fifth section, we define property A and prove an analogous
result for property A:

\begin{thm} Let $G(\sY)$ be a developable complex of groups over a
finite scwol $\sY$ with development $\sX$ such that
$\as\sX<\infty,$ and such that every local group $G_\sigma$ has
property A. If $\pi_1$ denotes the fundamental group of the
complex of groups, then $\pi_1$ has property A.
\end{thm}

In the final section, we use an example of a finitely generated
group which does not have either property A or finite $\as$ to
show that the results we obtain here cannot be improved by
relaxing the condition that the development have finite $\as.$

\section{Complexes of Groups}

Our notation and development will follow Bridson-Haefliger
\cite{BH}.

\begin{defn}  A {\em small category without loops} (abbreviated
{\em scwol}) is a set $\sX$ which is the disjoint union of a {\em
vertex set} $V(\sX)$ and an {\em edge set} $E(\sX).$  There are
maps \[i:E(\sX)\to V(\sX)\] and \[t:E(\sX)\to V(\sX)\] which
assign to each edge $a$ the initial vertex of $a$ and the terminal
vertex of $a,$ respectively.  Let $E^{(2)}(\sX)=\{(a,b) \in
E(\sX)\times E(\sX)\mid i(a)=t(b)\}$ denote the pairs of
composable edges. There is also a map \[E^{(2)}(\sX)\to E(\sX)\]
which assigns to each pair $(a,b)$ an edge $(ab)$ called the
composition of $a$ and $b.$  These maps are required to satisfy:
\begin{enumerate}
\item $i(ab)=i(b),$ and $t(ab)=t(a)$ for all $(a,b)\in
E^{(2)}(\sX);$
\item $a(bc)=(ab)c$ for all edges $a,b$ and $c$ with $i(a)=t(b)$
and $i(b)=t(c);$ and
\item $i(a)\neq t(a).$ (the no loops condition)
\end{enumerate}
\end{defn}

Let $E^{(k)}(\sX)$ denote the composable sequences of edges of
length $k,$ i.e., $E^{(k)}(\sX)=\{(a_1,\ldots,a_k)\in
(E(\sX))^k\mid i(a_i)=t(a_{i+1}),$ for $i=1,\ldots, k-1\}.$  By
convention, $E^{(0)}(\sX)=V(\sX).$ We define the {\em dimension}
of the scwol $\sX$ to be the maximum $k$ such that $E^{(k)}(\sX)$
is not empty.

\begin{defn}  The {\em geometric realization} $|\sX|$ is a piecewise
Euclidean polyhedral complex, with each $k$-cell isometric to the
standard simplex $\Delta^k.$  There is one such $k$-simplex $A$
for each $A\in E^{(k)}(\sX).$  The identifications are the obvious
ones, induced by the face relation among simplices.
\end{defn}

Observe that the geometric realization need not be a simplicial
complex, since it may be the case that the intersection of two
simplices is a union of faces.  One can eliminate this problem by
taking the barycentric subdivision, if one requires simplicial
complexes.  The geometric dimension of $|\sX|$ is the same as the
dimension of the combinatorial object $\sX.$

\begin{defn}  A {\em group action} on a scwol is a homomorphism $G\to
Aut(\sX)$ satisfying
\begin{enumerate}
\item for every $g\in G,$ and for all $a\in E(\sX)$ $g.i(a)\neq
t(a).$
\item for every $g\in G,$ and for all $a\in E(\sX)$ if $g.
i(a)= i(a),$ then $g. a=a.$
\end{enumerate}
\end{defn}

Notice that a group action on a scwol induces an isometric action
of the group on the geometric realization $|\sX|.$  Since we are
primarily concerned with isometric actions on metric spaces, this
is the action that we consider.

One forms the quotient $\sY=G\backslash \sX$ of the scwol $\sX$ by
the action of $G$ by taking $V(\sY)=G\backslash V(\sX),$ and
$E(\sY)= G\backslash E(\sX).$  One can verify that $\sY$ has the
structure of a scwol.

\begin{defn}  A {\em complex of groups $G(\sY)$ over a scwol} $\sY$ is a
collection $G(\sY)=(G_\sigma, \psi_a, g_{a,b})$ satisfying
\begin{enumerate}
\item to each $\sigma\in V(\sY),$ there corresponds a group
$G_\sigma$ called the {\em local group} at $\sigma;$
\item  for each $a\in E(\sY)$ there exists an injective
homomorphism $\psi_a:G_{i(a)}\to G_{t(a)};$ and
\item For each $(a,b)\in E^{(2)}(\sY),$ there is a $g_{a,b}\in
G_{t(a)}$ such that \newline\indent\indent(i)
$Ad(g_{a,b})\psi_{ab}=\psi_a\psi_b,$ where $Ad(g_{a,b})$ denotes
conjugation by $g_{a,b},$ and
\newline\indent\indent (ii)
$\psi_a(g_{b,c})g_{a,bc}=g_{a,b}g_{ab,c},$ for all $(a,b,c)\in
E^{(3)}(\sY).$
\end{enumerate}
\end{defn}

Given a group $G$ and an action of $G$ on the scwol $\sX,$ there
is an explicit construction of the complex of groups over the
quotient scwol which we do not describe here.  However, on the
other hand, it is not always the case that an arbitrary complex of
groups can be associated to a group action on some scwol $\sX.$
When this occurs, we say that the complex of groups is
\emph{developable} and we refer to the associated scwol $\sX$ as
the development of $G(\sY).$

It is clear that scwols of dimension 1 must have precisely two
types of vertices: sources and sinks.  A source is an initial
vertex of every edge it is contained in and a sink is a terminal
vertex of every edge it is contained in. Every one-dimensional
simplicial complex (graph) can be given the structure of a
one-dimensional scwol by placing a source vertex in the middle of
every edge, thus giving the original vertices the structure of
sinks.  It is easy to verify that the theory of complexes of
groups over one-dimensional scwols is precisely the same as the
theory of graphs of groups. Phrased in terms of the language of
complexes of groups, the Bass-Serre structure theorem for groups
acting without inversion on graphs says that if $dim(\sY)=1,$ then
$G(\sY)$ is always developable.

When a complex of groups is developable, there is an explicit
method of constructing both the scwol $\sX$ and the group $G$
which acts on the scwol.  The scwol $\sX$ on which the group acts
is simply connected and has an explicit description in a similar
way to the construction of the tree $\tilde X$ in the theory of
graphs of groups (see \cite{S}).

Indeed, if $G(\sY)$ is a developable complex of groups, then we
can define the development $D(\sY)$ to be the scwol whose vertices
and edges are given by $V(D(\sY))=\{(gG_\sigma, \sigma)\mid
\sigma\in V(\sY)\},$ and $E(D(\sY))=\{(gG_{i(a)},a)\mid a\in
E(\sY)\}.$ Then the group $G$ acts on the development $D(\sY)$ by
left multiplication. The development is isomorphic to the scwol
$\sX,$ mentioned above.  (See \cite{BH} for more details.)

We describe the fundamental group of the complex of groups
$\pi_1(G(\sY))$ which is the group $G,$ up to isomorphism.  As in
the theory of graphs of groups, there are two equivalent
descriptions of the fundamental group.  Both rely on the
construction of the auxiliary group $FG(\sY).$ Let $E^\pm(\sY)$
denote the collection of symbols $\{a^+,a^-\}$ where $a\in
E(\sY).$  The elements of $E^\pm(\sY)$ can be thought of as {\em
oriented edges}.  If $e=a^+,$ then define $i(e)=t(a)$ and
$t(e)=i(a).$  Accordingly, if $e=a^-,$ define $t(e)=t(a)$ and
$i(e)=i(a).$ Then, define $FG(\sY)$ to be the free product of the
local groups $G_\sigma$ and the free group generated by the
collection $E^\pm(\sY)$ subject to the additional relations:
\begin{enumerate}
\item $(a^+)^{-1}=a^-,$ and $(a^-)^{-1}=a^+;$
\item $a^+b^+=g_{a,b}(ab)^+;$
\item $\psi_a(g)=a^+ga^-,$ for all $g\in G_{i(a)}.$
\end{enumerate}

The first description of the fundamental group is in terms of
$G(\sY)$-loops based at some fixed vertex $\sigma_0.$   An edge
path in $\sY$ is a sequence $(e_1,\ldots,e_k)$ with
$t(e_i)=i(e_{i+1}),$ for all $i=1,\ldots,k-1.$ By a $G(\sY)$-path
issuing from $\sigma_0$ we mean a sequence
$(g_0,e_1,g_1,\ldots,e_k,g_k),$ where $(e_1,\ldots,e_k)$ is an
edge path in $\sY,$ and $g_0\in G_{\sigma_0},$ and $g_i\in
G_{t(e_i)},$ for $i>0.$  We associate the word $g_0e_1\ldots
e_kg_k\in FG(\sY)$ to the path described above. A $G(\sY)$-loop
based at $\sigma_0$ is a $G(\sY)$ path with $t(e_k)=\sigma_0.$
There is an obvious notion of homotopy on $G(\sY)$-paths, i.e.,
the notion of homotopy on the geometric realization.  The
fundamental group $\pi_1(G(\sY),\sigma_0)$ is the collection of
all words associated to $G(\sY)$-loops based at $\sigma_0,$ up to
homotopy equivalence.

The second description is much simpler.  Let $T$ be a maximal tree
in $|\sY|^{(1)}.$  Then, $\pi_1(G(\sY),T)$ is $FG(\sY)$ subject to
the additional relation $a^+=1,$ for all $a\in T.$  For a
connected scwol, there is an isomorphism $\pi_1(G(\sY),\sigma_0)
\to \pi_1(G(\sY),T).$

\section{Groups Acting on Metric Spaces}

For group actions considered on a local scale, the stabilizer
plays a key r\^{o}le.  The corresponding notion for group actions
considered in the global sense is that of the $R$-stabilizer which
we define presently.

\begin{defn}  Let $\Gamma$ be a group acting on the pointed metric space
$(X, x_0)$ by isometries.  For every $R>0$ define the {\em
R-stabilizer} of the point $x_0,$ denoted $W_R(x_0),$ by
\[W_R(x_0)=\{\gamma\in\Gamma\mid d(\gamma x_0,x_0)\le R\}.\]
\end{defn}

In \cite{BD1}, Dranishnikov and the author characterized the
$R$-stabilizers of the action of a fundamental group of a graph of
groups on the tree corresponding to its development.  The
following proposition is a natural generalization of that result.

\begin{Proposition} \label{stab} Let $G(\sY)$ be a developable complex of
groups.  Fix a vertex $\sigma_0\in \sY,$ and consider the action
of $\pi_1(G(\sY),\sigma_0)$ on the simply connected scwol $\sX$
induced by the complex of groups.  Then the $R$-stabilizer
$W_R(\sigma_0)$ is the set of all elements of
$\pi_1(G(\sY),\sigma_0)$ with associated path $c$ of length not
exceeding $R.$
\end{Proposition}

\begin{proof}  Let $\gamma$ be a word in $\pi$ which is reduced
and for which the path length is $R.$  Then consider $d_\sX(\gamma
G_{\sigma_0},G_{\sigma_0}).$  The path corresponding to $\gamma$
lifts to a path in $\sX,$ so this distance is at most $R.$

Thus, it remains to show that no element of $\gamma G_{\sigma_0}$
has length less than $R$ if $\gamma$ has length $R$ in its most
reduced form.

To this end, let $xg=g_0e_1\cdots e_Rg_Rg,$ be a reduced word in
$xG_{\sigma_0}.$

Two of the relations on $FG(\sY)$ can affect the path length of a
word.  In order for the relation
$(e_i^+e_{i+1}^+)=g_{e_i,e_{i+1}}(e_ie_{i+1})^+$ to occur in this
word, we would need two composable edges, in the sense that
$t(e_{i+1})=i(e_i),$ but since by definition, we have
$t(e_i)=i(e_{i+1})$ for all edges, we conclude that
$e_{i}^{-1}=e_{i+1}.$  Thus, the only relation that can occur is
the type $\psi_a(g)=a^+ga^-.$

Suppose that a sequence of this type of relation occurs which
transform the word $xg$ into $g_0e_1\cdots g_{k-1}e_kg_kg^\prime
e_{k+1} g^\prime_k\cdots g_R^\prime.$
%In order for a type-(3) reduction
%to occur here, the edges $e_k$ and $e_{k+1}$ must be composable in
%the sense that $i(e_k)=t(e_{k+1}).$  But, by the definition of a
%path, $t(e_k)=i(e_{k+1}),$ so $e_k$ and $e_{k+1}$ are the same
%edge with opposite orientation.  Thus, the only possibility is a
%reduction of type (4).
%
Here two cases can occur.  In the first case, edge $e_k$ has
positive orientation.  Then, $e_{k}g_kg^\prime
e_{k+1}=\psi_{e_k}(g_kg^\prime).$  But, in this case, we can
obtain a reduction of the original path $x$ since
$e_kge_{k+1}=\psi_{e_k}(g_k).$  As the original path was reduced,
this cannot occur.  The other case is when $e_{k+1}$ has positive
orientation.  Then, $g_kg^\prime=\psi_{e_{k+1}}(h)$ for some $h.$
But, it must be the case that $g^\prime$ was generated by this
type of relation, as this is the only type that can occur and
involves edges.  Thus, $g^\prime$ is itself
$\psi_{e_{k+1}}(h^\prime),$ for some $h^\prime.$  Thus,
$g=\psi_{e_{k+1}}(h (h^\prime)^{-1}).$ This enables a reduction of
the original word.  Thus, the length of any representative of
$xG_{\sigma_0}$ has length no less than $R$ as desired.
\end{proof}

\section{Asymptotic Dimension}

Asymptotic dimension was introduced by Gromov \cite{Gr1}. It is
the coarse analog of Ostrand's characterization of covering
dimension for metric spaces, \cite{Ost}.

\begin{defn}  Let $X$ be a metric space.  We define the {\em asymptotic
dimension} $(\as)$ of $X$ by the following inequality: $\as X\le
n$ if for every $R>0$ there exist $n+1$ families of sets
$\sU_0,\ldots,\sU_n$ which are uniformly bounded, which cover $X$
and which are $R$-disjoint in the sense that any two distinct sets
from the same family are at a distance greater than $R$ from each
other.  We define $\as X=n$ if it is the case that $\as X\le n,$
but it is not the case that $\as X\le n-1.$
\end{defn}

%The finiteness of $\as$ for a finitely generated group became
%important following the paper of Yu \cite{Yu1} where it is shown
%that the Novikov higher signature conjecture holds for groups with
%finite $\as.$  Originally it was not known whether there were
%finitely generated groups with infinite $\as,$ but recent examples
%of Gromov and others \cite{Gr2}, \cite{Detal} have made the
%question of determining precisely which groups have finite $\as$
%important.  See the final section for further consideration of
%this group.

The goal of this section is to see that the finiteness of $\as$ is
preserved by the construction of the fundamental group of a
developable complex of groups.  This is a natural generalization
of the main theorem in \cite{BD2}.

\begin{defn}  Let $X_\alpha$ be a family of subsets of the metric space
$X.$ We say that the family satisfies the inequality $\as
X_\alpha\le n$ {\em uniformly} if for every $D>0$ there is a
number $R>0$ and a collection of $D$-bounded, $R$-disjoint
families $\{\sU_i^\alpha\}$ so that for each $\alpha,$
$\{\sU_i^\alpha\}$ covers $X_\alpha.$
\end{defn}

A common example of a family satisfying $\as X_\alpha\le n$
uniformly is a family of isometric metric spaces.

The following union theorem appears as Theorem 1 in \cite{BD1}.

\begin{thm}[Union Theorem] Let $X=\cup_\alpha F_\alpha$ and
$\as F_\alpha$ uniformly.  Suppose that for any $r$ there exists a
set $Y_r\subset X$ such that $\as Y_r\le n$ and the family
$\{F_\alpha\setminus Y_r\}$ is $r$-disjoint.  Then, $\as X\le n.$
\end{thm}

As a corollary, we have the following finite union theorem.

\begin{thm}[Finite Union Theorem]  Let $X=\cup_{i=1}^k X_i$ be a
metric space.  Then, $\as X\le \max\{\as X_i\mid i=1,\ldots k\}$
\end{thm}

Let $G(\sY)$ be a developable complex of groups, with $\sY$
finite and $\as |\sX|\le k.$  Suppose further that the local
groups are finitely generated. Then, the finiteness of $\sY$
implies that the fundamental group $\pi$ is finitely generated.
Indeed, if $S_\sigma$ denotes a finite generating set for each
local group $G_\sigma,$ then we can consider $FG(\sY)$ in the
metric obtained from the disjoint union of all the $S_\sigma$ and
the set $E^\pm(\sY).$  Thus, the notion of $\as$ is well-defined
for the fundamental group of a complex of groups.

\begin{Lemma} \label{asdim} Let $\pi$ be the fundamental group of a complex of
groups $G(\sY)$ where $\sY$ is finite, the local groups are
finitely generated, and the local groups satisfy $\as G_\sigma\le
n,$ then for every $R>0,$ $\as W_R(G_{\sigma_0})\le n.$
\end{Lemma}

\begin{proof}  In Proposition 1 we characterized $W_R(G_{\sigma_0})$ as
the set of all elements in $\pi$ with length at most $R.$  Let
$\sX$ denote the development of $FG(\sY).$

In order to apply an inductive argument, we consider a larger set,
$K\subset FG(\sY)$ which is the set of all words in $FG(\sY)$
issuing from $\sigma_0.$  The group $\pi$ is a subset of $K$ and
the set $K$ acts on $\sX$ by left multiplication.  We show that
the $R$-stabilizer of this action has $\as$ at most $n.$  It
follows then that the $R$-stabilizer of the action of $\pi$ on
$\sX$ will also have $\as$ at most $n.$

In light of the finite union theorem, in order to show $\as
W_R(G_{\sigma_0})\le n,$ it suffices to show that the subset
$K_j\subset K$ of reduced words in $K$ with length equal to $j$
has $\as K_j\le n.$ Indeed, $W_R(G_{\sigma_0})\subset \cup_{j=0}^R
K_j,$ which is a finite union.

We proceed by induction.  The case $j=0$ is true by assumption
since $K_0=G_{\sigma_0}.$  Consider the case $K_{j+1}$ with $j\ge
0.$  Observe that $K_{j+1}\subset \cup_{a\in E^\pm(\sY)} K_j
aG_{t(a)}.$

The orientation of the edge $a$ is an issue, as it determines
whether the group $G_{t(a)}$ is a domain or codomain of the
function $\psi_a.$  Thus, it is necessary to consider two cases
separately.

Suppose first that $a$ has negative orientation.  So, we are
considering $K_ja^-G_{t(a)}.$  For every $r>0$ let
$Y_r=K_ja^-N_r(\psi_a(G_{i(a)})),$ where the $r$-neighborhood is
taken in the group $FG(\sY).$  Then $Y_r$ is coarsely equivalent
to $K_ja^-\psi_a(G_{i(a)}).$  Now we have
$K_ja^-\psi_a(G_{i(a)})=K_jG_{i(a)} a^-,$ which is just $K_ja^-.$
Finally, as $K_ja^-$ is coarsely equivalent to $K_j,$ we have $\as
Y_r=\as K_j,$ which by the inductive hypothesis does not exceed
$n.$

Next, decompose the set $K_ja^-G_{t(a)}$ into families
$\{xa^-G_{t(a)}\},$ where the index runs over all $x\in K_j$ which
do not end with an element $g\in G_{i(a)}.$  One can still obtain
these elements through the relations of $FG(\sY).$  For instance,
to obtain $xga^-g^\prime,$ with $x$ of the required form, $g\in
G_{i(a)}$ and $g^\prime$ in $G_{t(a)},$ simply take the word
$xa^-\psi_a(g)g^\prime,$ which is of the required form.  Next,
observe that the map $G_{t(a)}\mapsto xa^-G_{t(a)}$ is an isometry
in the (left-invariant) word metric, so the family
$\{xa^-G_{t(a)}\}$ has $\as \le n,$ uniformly.

In order to apply the union theorem to this family, it remains to
show only that the family $\{xa^-G_{t(a)}\setminus Y_r\}$ is
$r$-disjoint.  To this end, let $xa^-z$ and $x'a^-z'$ be given in
different families.  Then we compute
$d(xa^-z,x'a^-z')=\|z^{-1}a^+x^{-1}x'a^-z'\|.$ Since $z$ and $z'$
lie outside of $N_r(\psi_a(G_{i(a)})),$ take $z=s\psi_a(g)$ and
$z'=s'\psi_a(g'),$ where $\|s\|>r,$ $\|s'\|>r,$ and $\psi_a(g)$
and $\psi_a(g')$ are in $\psi_a(G_{t(a)}).$  Then,
\[\|z^{-1}a^+x^{-1}x'a^-z'\|=\|s^{-1}a^+g^{-1}x^{-1}x'g'a^-z'\|.\]
Now, in order for this length to be less than $r,$ a reduction
must occur in the middle, so that $a^+$ and $a^-$ annihilate each
other.  In order for this to occur, we must have
$g^{-1}x^{-1}x'g'\in G_{i(a)}.$  Thus, $x^{-1}x'\in G_{i(a)}.$
But, this means that $xa^-G_{t(a)}$ and $x'a^-G_{t(a)}$ define the
same family.  Thus, in the case that the edge has negative
orientation, we have $\as K_ja^-G_{t(a)}\le n.$

Next, we consider the case where the edge $a$ has positive
orientation.  In this case,
$K_ja^+G_{i(a)}=K_j\psi_a(G_{i(a)})a^+,$ which is coarsely
equivalent to $K_j.$  We conclude that $\as K_ja^+G_{i(a)}=\as
K_j\le n.$
\end{proof}

The following result appears as Theorem 2 from \cite{BD1}.

\begin{thm} Assume that a finitely generated group
$\Gamma$ acts by isometries on a metric space $X$ with a base
point $x_0$ and with $\as X\le k.$  Suppose that $\as W_R(x_0)\le
n$ for all R.  Then $\as\Gamma\le (n+1)(k+1)-1.$
\end{thm}

This estimate on the dimension is far from sharp.  It is useful
only as a means to prove that $\as\Gamma<\infty.$  The exact
estimate should be $n+k.$ (See \cite{BD2} for the proof of the
exact formula in the case of groups acting on trees by
isometries.)

As a consequence of the preceding theorem, we have our main result
on $\as.$

\begin{Theorem}
Let $\Gamma$ be the fundamental group of a finite developable
complex of groups $G(\sY)$ corresponding to an action by
isometries on the geometric realization of the scwol $\sX.$
Suppose that the local groups are finitely generated and that $\as
G_\sigma\le n.$ Assume additionally that $\as |\sX|\le k.$ Then
$\as \Gamma\le (n+1)(k+1)-1.$
\end{Theorem}

The important result is summarized as a corollary:

\begin{cor}  Let $\Gamma$ be the fundamental group of a
finite developable complex of groups $G(\sY)$ such that the
development $\sX$ has $\as|\sX|<\infty,$ and such that every base
group $G_\sigma$ has $\as G_\sigma<\infty.$  Then,
$\as\Gamma<\infty.$
\end{cor}

\section{Property A}

Property A was introduced by G. Yu, \cite{Yu2}.  It is a weak form
of amenability which, for groups, implies the existence of a
uniform embedding into Hilbert space.  Thus, the coarse
Baum-Connes conjecture and the Novikov conjecture hold for this
group.

\begin{defn}  Let $X$ be a metric space.  Let $P(X)$ denote the
set of probability measures on $X$ in the $l_1$ metric.  The
metric space $X$ has property $A$ if there exists a sequence of
maps $a^n:X\to P(X)$ satisfying the following two conditions:
\begin{enumerate}
\item for every $n$ there is an $R$ so that for every $x,$
$\supp a^n_x\subset B_R(x)$, and
\item for every $K>0,$
\[\lim_{n\to\infty}\sup_{d(x,y)<K}\|a^n_x-a^n_y\|_1=0.\]
\end{enumerate}
\end{defn}

As mentioned in the introduction finitely generated groups with
finite $\as$ have property A, \cite{HR}.  Tu proved, \cite{Tu},
that the fundamental group of a finite graph of groups in which
each vertex group has property A will have property A. In
\cite{B}, the author generalized Tu's results to groups acting by
isometries on metric spaces with finite $\as.$

In particular, the theorem proved in \cite{B} is the following:

\begin{thm} Assume that the finitely generated group $\Gamma$ acts on
the metric space $X$ by isometries.  Assume that $\as X\le n,$ and
that for every $R,$ the $R$-stabilizer of a basepoint $x_0\in X$
has property A.  Then $\Gamma$ has property A.
\end{thm}

In \cite{B} the author proves a union theorem which is analogous
to the union theorem for $\as$ from $\S 4.$

\begin{thm}[Union Theorem]  Let $X=\cup_\alpha X_\alpha$ where the
$X_\alpha$ are pairwise isometric and have property A.  Suppose
further that for every $r>0$ there is a set $Y_r$ so that
$\{X_\alpha \setminus Y_r\}$ is $r$-disjoint.  Then, $X$ has
property A.
\end{thm}

As a consequence we obtain the finite union theorem.

\begin{thm}[Finite Union Theorem]  Let $X=\cup_{i=1}^n X_i$ where the $X_i$ all have
property A.  Then, $X$ has property A.
\end{thm}

The analog of Lemma \ref{asdim} for Property A is the following.
The proof follows from the proof of Lemma \ref{asdim}.

\begin{Lemma}  Let $G(\sY)$ be a developable complex of groups. Suppose
that $\sY$ a finite, connected scwol, that the local groups are
finitely generated, and that the local groups have property A.
Then, the $R$-stabilizer $W_R(\sigma_0)$ for some base vertex
$\sigma_0$ also has property A.
\end{Lemma}

Applying the theorem from \cite{B} cited above, on groups acting
on metric spaces with finite $\as,$ we obtain the following
generalization of Tu's theorem.

\begin{Theorem}  Let $G(\sY)$ be a developable complex of groups with
corresponding development $\sX$ and fundamental group $\pi.$
Suppose that $\as |\sX|$ is finite and that the stabilizers of the
action have property A.  Then, $\pi$ has property A.
\end{Theorem}

\section{Example}

The following example illustrates that one must consider the
large-scale structure of the development.

Consider a finitely presented group $\Gamma$ which does not have
finite $\as$ or property A, see \cite{Gr2}, \cite{Detal}. Since
the group is finitely presented, there is a finite complex $K$ so
that $\pi_1(K)=\Gamma.$  Thus, by taking the complex of groups
with each vertex trivial and the scwol whose geometric realization
is equal to the complex $K,$ one obtains the fundamental group of
the complex of groups equal to the group $\Gamma.$  The vertex
groups have finite $\as$ and the group acts on the universal
cover, so the complex of groups is developable. The complex $K$ is
finite, yet the group $\Gamma$ does not have finite $\as$ and does
not have property A.

\end{document}